\documentclass[letterpaper, 10pt, conference]{ieeeconf}

\IEEEoverridecommandlockouts                              % This command is only
\let\labelindent\relax

\usepackage[cmex10]{amsmath}
\usepackage{mathrsfs}
\usepackage{empheq}
\usepackage{tikz}
\usetikzlibrary{fadings}
\tikzfading[name=fade out,
            inner color=transparent!20,
            outer color=transparent!90]
\tikzfading[name=fade equal,
            inner color=transparent!30,
            outer color=transparent!30]

\usepackage{amssymb}
\usepackage{txfonts}
\let\labelindent\relax
\usepackage{enumitem}
\usepackage{dsfont,citesort}
\usepackage{tikz}
\usetikzlibrary{automata, positioning}
\usepackage{subfig}
\usepackage{mathtools}
\allowdisplaybreaks
\newtheorem{theorem}{Theorem}
\newtheorem{corollary}[theorem]{Corollary}
\newtheorem{lemma}[theorem]{Lemma}
\newtheorem{proposition}[theorem]{Proposition}
\newtheorem{definition}[theorem]{Definition}
\newcommand{\cE}{{\mathcal E}}
\newcommand{\cK}{{\mathcal K}}
\newcommand{\cS}{{\mathcal S}}

\newcommand{\ree}{\mathbb{R}}
\newcommand{\im}{\textup{Im}}

\newcommand{\E}{{\cal E}}

\allowdisplaybreaks

\usepackage{url}

\newcommand{\R}{\ensuremath{{\mathbb R}}}
\newcommand{\N}{\ensuremath{{\mathbb N}}}

\DeclareMathOperator{\argmin}{argmin}

\colorlet{Darkred}{red!50!black}
\colorlet{Darkgreen}{green!50!black}
\xdefinecolor{UIUCblue}{rgb}{0,0.24,0.49} %professional:{110,139,191}
\xdefinecolor{UIUCorange}{rgb}{0.96,0.5,0.14} %professional:{239,138,28}

\newlist{ass}{enumerate}{1}
\setlist[ass]{label = {\bf (A\arabic*)}, resume}

\title{ \LARGE \bf
	Existence and Completeness of Solutions to  Extended Projected Dynamical Systems and Sector-bounded Projection-based  Controllers 
}

\author{ W.P.M.H. Heemels \and Aneel Tanwani  %
\thanks{This research is supported by the ERC Advanced Grant PROACTHIS.}
\thanks{Maurice Heemels is at Control Systems Technology Section, Dept.~Mechanical Eng., Eindhoven Univ.~Technology. Email: \texttt{m.heemels@tue.nl}}
\thanks{Aneel Tanwani is with CNRS -- LAAS, Toulouse. }%Email: \texttt{aneel.tanwani@cnrs.fr}}
}

\begin{document}

\maketitle
\begin{abstract}
Projection-based control (PBC) systems have significant engineering impact and receive considerable scientific attention. To properly describe closed-loop PBC systems, extensions of classical projected dynamical systems are needed, because partial projection operators and irregular constraint sets (sectors) are crucial in PBC. These two features obstruct the application of existing results on existence and completeness of solutions. To establish a rigorous foundation for the analysis and design of PBC, we provide essential existence and completeness properties for this new class of discontinuous  systems. 
\end{abstract}

\begin{keywords}
Projected dynamical systems, hybrid control,  discontinuous dynamics, sectors,  non-smooth analysis. 
\end{keywords}

\thispagestyle{plain}
\pagestyle{plain}

\section{Introduction}

In recent years, there has been a strong interest in projection-based control (PBC) systems, including the { hybrid integrator-gain system} (HIGS) \cite{DeeSha21,EijHee_CSL20a,ShiNi_CDC22}, in which specific closed-loop signals are kept in sector-bounded sets in order to overcome fundamental performance limitations inherent to LTI control \cite{Seron1997}. PBC's potential in overcoming these limitations was demonstrated in \cite{EijHee_CSL20a} and practical successes were reported in lithography systems \cite{Deenen2017}, force atomic microscopes \cite{ShiNi_CDC22}, etc. Interestingly, the fundamental and important problem of well-posedness of PBC in the sense of existence and completeness of  solutions was so far only partially addressed in \cite{ShaHee_CDC19a,DeeSha21}, where the plant model was limited to be LTI and the external input signals (references, disturbances) were restricted to so-called (piecewise) Bohl functions (i.e., functions generated by LTI models). Clearly, to expand the applicability of PBC, there is a strong interest to provide well-posedness guarantees for larger classes of nonlinear plants and controllers, and more general and natural input classes, such as piecewise continuous inputs. 

The natural formalism to describe (closed-loop) PBC systems is formed by the recently proposed extension \cite{ShaHee_CDC19a,DeeSha21} of classical projected dynamical systems (PDS) \cite{dupuis1993dynamical,nagurney2012projected}.  ``Classical'' PDS consider a differential equation given by 
\begin{equation} \label{eq:f}
\dot x(t) = f(x(t))
\end{equation}
in which the state $x(t)\in \ree^n$ is restricted to remain inside a set $S\subseteq \ree^n$, for all times $t\in\ree_{\geq 0}$, which in PDS is ensured by redirecting the vector field at the boundary of $S$. Formally, PDS are given for a continuous vector field $f: \ree^n \rightarrow \ree^n $ and set $S\subseteq \ree^n$ (with further additional conditions) 
by 
\begin{align} \label{pds1}
\dot x &= \Pi_{S}(x,f(x)) \text{ with} \\ 
\Pi_{S}(x,v) &= \argmin_{w \in T_S(x)} \|w-v\|  \label{PiS0}
\end{align}
for $x\in S$ and $v\in \ree^n$. The set $T_S(x)$ is the tangent cone of $S$ at $x$, defined formally in Section~\ref{sec:epds}, which essentially contains all admissible velocities that keep the trajectories inside $S$. It is possible to draw connections between this conventional PDS and other formalisms used for modeling constraints on evolution of state trajectories, which provide some insight about the algorithms for simulating such systems \cite{brogliatoSCL}. %Control design problems using such formalisms, with the objective of constraining the motion of the state to a prescribed set have been studied in \cite{TanwBrog18}. %Using similar tools, one can write analogues for \eqref{piSE} where the geometry of the subspace $\cE$ also plays a role. 

Although the framework of PDS is a source of inspiration to study PBC, it cannot properly describe the resulting closed-loop PBC systems. This can be observed from the fact that \eqref{PiS0} allows projection  along all possible directions of the state vector (including both controller and plant states) in the sense that it just takes the vector 
$\Pi_{S}(x,v) \in T_S(x)$ that is ``closest" to $v$ irrespective of the direction $\Pi_{S}(x,v)-v$.  Clearly, if  \eqref{eq:f} is a closed-loop system in the sense of an interconnection of a physical plant and a controller (and thus the state $x$ consists of physical plant states $x_p$ and controller states $x_c$), one cannot project freely in all directions. Indeed, the physical state dynamics cannot be modified by projection; it is only possible to ``project" the controller ($x_c$-)dynamics.  Hence, in contrast to PDS, there are only limited directions in order to ``correct'' the vector field $f(x)$ at the boundary, which we describe by a subspace $\cE \subseteq \ree^{n}$. One then obtains
\begin{align} 
\dot x & = \Pi_{S,\cE}(x,f(x)) \text{ with }\label{pds2} \\ 
\Pi_{S,\cE}(x,v) &= \argmin_{w \in T_S(x), w - v \in \cE} \|w-v\| \label{piSE}
\end{align}
for $x\in S$ and $v\in \ree^n$. Hence, the projection $\Pi_{S,\cE}$ projects the vector $v$ onto the set of admissible velocities (in the tangent cone $T_S(x)$, $x \in S$) along $\cE$ in such a way that the correction $w-v$ is minimal in norm. For these systems, we coined the term {\em extended Projected Dynamical Systems} (ePDS) in \cite{ShaHee_CDC19a,DeeSha21}, as they include the classical PDS \eqref{pds1} as a special case by taking $\cE=\ree^n$. Clearly, this new projection operator \eqref{piSE} and the corresponding ePDS \eqref{pds2} require careful analysis to provide conditions on the vector field $f$, constraint set $S$ and projection directions $\cE$ so that $\Pi_{S,\cE}$ is well-defined, and the existence and completeness of solutions to \eqref{pds2} is guaranteed.

Recently, such questions were answered for an alternative extension \cite{hauswirth2021projected} of PDS, called oblique PDS, that did not restrict the projection directions as is needed for closed-loop PBC systems, but allowed a state-dependent metric to execute the projection in \eqref{PiS} leading to, loosely speaking, 
\begin{equation} \label{PiS}
\Pi_{S}(x,v) = \argmin_{w \in T_S(x)} \|w-v\|_{G(x)}, 
\end{equation}
where $G(x)$ is a positive definite matrix for each $x\in S$ and $\|w-v\|_{G(x)}^2 := (w-v)^\top {G(x)} (w-v)$. Although under rather strict conditions \cite{ShaHee_ACC21a} some connections can be established between extended and oblique PDS, the underlying philosophy is different, as well as the underlying mathematical structure. Moreover, the well-posedness conditions in \cite{hauswirth2021projected} are given for constraint sets that are convex or Clarke regular \cite{rockafellar2}, while in the PBC, including HIGS \cite{DeeSha21,EijHee_CSL20a,ShiNi_CDC22,Deenen2017,ShiNi_CDC22},  sectors are used as constraint sets that do not satisfy these regularity properties. This calls for new and alternative conditions guaranteeing the well-posedness of the ePDS \eqref{pds2} with partial projection operators as in \eqref{piSE} and irregular constraint sets such as sectors. Such results will form the main contributions of the present paper, thereby laying a rigorous foundation for the analysis and design of PBC. The first set of well-posedness results will be given in Section~\ref{sec:epds} for ePDS, hence, with partial projection, but still working with regular sets being finitely generated. Based on these results, in Sections~\ref{sec:epds_sectors} and~\ref{sec:well}, we will provide well-posedness results for ePDS with sectorsets and partial projection, and will also show how closed-loop PBC systems are covered by these results.

\section{Extended PDS on Finitely Generated Constraint Sets} \label{sec:epds}

The primary object in our study of ePDS is the operator $\Pi_{S,\cE}(x,f(x))$, which basically projects the unconstrained vector field $f(x)$ on the set $T_S(x)$, in the direction determined by $\cE$, for each $x \in S$. The tangent cone to a set $S \subset \mathbb{R}^n$ at a point $x \in S$, denoted by $T_S(x)$, is the set of all vectors $v\in \mathbb{R}^n$ for which there exist sequences $\{x_i\}_{i\in \mathbb{N}} \in S$ and $\{\tau_i\}_{i\in \mathbb{N}}$, $\tau_i > 0$
with $x_i \rightarrow x$, $\tau_i \downarrow 0$ and $i \rightarrow \infty$, such that $
    v = \lim_{i\rightarrow \infty} \frac{x_i - x}{\tau_i}.$

For ease of exposition, we consider  sets $S$, which are described by the intersections of sublevel sets of finitely many real-valued functions. In particular, for given functions $h_i:\R^n \to \R$, $i=1,\dots,m$, we take $S$ to be of the form
\begin{equation}\label{eq:defSFinite}
S:= \{ x\in \R^n \, \vert \, h_i(x) \ge 0, \text{ for all } i=1,\dots,m\}.
\end{equation}
For all $x \in \R^n$, we define the set of active constraints by
\begin{equation} \label{eq:activeconstraints}
    J(x) = \{ i \in \{1,\dots,m\} \mid  h_i(x) = 0 \}.
\end{equation}
{\bf (CQ)} The functions $h_i$, $i=1,\dots,m$, are assumed to be smooth and for each $x\in S$, it holds that $\{\nabla h_i(x), i \in J(x)\}$ are linearly independent. %$\nabla h_i$ not vanishing in a neighborhood of $\{x \in \R^n \, \vert \, h_i(x) = 0\}$. %{\bf MH: only elementwise non-vanishing is OK, right? no linearly independence assumptions on $\nabla h_{J(x)}$ is needed, right?}

Under this constraint qualification, the tangent cone to $S$ at $x$ is given by \cite[Lemma~12.2]{NocWri00}, \cite{brogliatoSCL}:
\begin{equation} \label{eq:Tsexplicit}
    T_S(x) = \{v \in \R^n \, \vert \, \langle \nabla h_i(x), v \rangle \ge 0, i \in J(x)\}.
\end{equation}

%\begin{remark}
%There are different notions of tangent cones in the literature. The one considered is also called Bouligand's tangent cone and for our purposes, this definition suffices. The set $S$ defined in \eqref{eq:defSFinite} has the property that several notions of the tangent cones from the literature coincide, and in particular, they are Clarke regular  and the tangent cone $T_S(x)$ is convex for each $x \in S$.
%\end{remark}
{\bf Well-defined projection operator $\Pi_{S,\cE}$:} For the dynamical system  \eqref{pds2}, it is important to study conditions under which the right-hand side is well-defined, i.e.,  $\Pi_{S,\cE}(x,f(x))$ is  nonempty, and preferably single-valued, for each $x \in S$ and $f(x) \in \R^n$. %, although it can be multi-valued. 
Clearly, for a given $S$ and $\cE$, this may not be the case in general, so  we need suitable conditions. %$\Pi_{S,\cE}(x,f(x))$ is nonempty for a given subspace $\cE$ and $S$ described in \eqref{eq:defSFinite}.

\begin{proposition}
    Consider a closed set $S$, and a given subspace $\cE \subseteq \R^n$. For each $x\in S$, if
    \begin{equation} \label{eq:nonempty}
    T_S(x) \cap (f(x) + \cE ) \neq \emptyset
    \end{equation}
    then $\Pi_{S,\cE}(x,f(x))$ is non-empty. If $S$ satisfies \eqref{eq:defSFinite} and {\bf(CQ)} holds, then \eqref{eq:nonempty} implies that $\Pi_{S,\cE}(x,f(x))$ is a singleton.
\end{proposition}

{\bf Existence of Solutions: } We now turn our attention to the existence of Carath\'eodory solutions to the ePDS \eqref{pds2}. 

\begin{definition}
We call a function $x:[0,T]\rightarrow \ree^{n}$ a (Carath\'eodory) solution to  \eqref{piSE}, if $x$ is absolutely continuous on $[0,T]$ and satisfies $\dot x(t) = \Pi_{\cS,\cE}(x(t),f(x(t)))$ for almost all $t\in[0,T]$ and $x(t)\in \cS$ for all $t\in[0,T]$. %{\bf should we add the latter or is it implied?}. 
We say that $x:[0,\infty)\rightarrow \ree^{n}$ is a solution on $[0,\infty)$, if the restriction of $x$ to $[0,T]$ is a solution on $[0,T]$ for each $T>0$.
\end{definition}

Our road map  to establish existence of Carath\'eodory solutions  is based on constructing the Krasovskii regularization of the {\em discontinuous} dynamical system \eqref{pds2} and demonstrate the existence of solutions to this regularization with additional viability conditions in the sense of Aubin \cite{aubin}. We show then that all Krasovskii solutions satisfy the viability condition, see \eqref{eq:viab=PDS} in the forthcoming Theorem~\ref{thm.kras=pds}, which essentially says that the Carath\'eodory solutions coincide with Krasovskii solutions for the ePDS. 

To state the result, we let $F(x):= \Pi_{S,\cE}(x,f(x))$, and denote the Krasovskii regularization by $K_F(x):= \cap_{\delta >0} \overline{\rm con}F(B(x,\delta))$, where $\overline{\textup{con}}(M)$ denotes the closed convex hull of the set $M$, in other words, the smallest closed convex set containing $M$. 

\begin{theorem} \label{thm.kras=pds}
    Assume that $f$ is continuous and the set $S$ in \eqref{eq:defSFinite} satisfies ${\bf(CQ)}$ and \eqref{eq:nonempty}.
    Then, for all $x\in S$, it holds that 
    \begin{equation}\label{eq:viab=PDS}
        K_F(x) \cap T_S(x) = \{F(x)\} =  \{\Pi_{S,{\cE}}(x,f(x)) \}.
    \end{equation}
\end{theorem}
For the proof of this result, we need to compute the Krasovskii regularization of $F$, which is described as follows:

\begin{proposition} \label{prop:kras}
For a continuous function $f:\R^n \to \R^n$, a closed set $S \subset \R^n$ satisfying \eqref{eq:nonempty}, and $F(x)= \Pi_{S,\cE}(x,f(x))$, it holds that 
\begin{equation}\label{eq:convKrasovskii}
K_F(x) = {\rm con } \limsup_{y \to x} P_{T_S(y),\cE} f(y),
%K_F(x) = {\rm con} \limsup_{y \to x} F(y)
\end{equation}
where $P_{K,\cE} (f(y)) := \argmin_{v \in K, f(y) - v \in \cE} \| v - f(y) \|$,  and the lim~sup on the right is interpreted in terms of convergence of sets. 
In particular, if $S$ satisfies \eqref{eq:defSFinite} and {\bf(CQ)} holds, then 
\begin{equation}\label{eq:krasovskiiFinite}
    K_F(x) = \textup{con}\{P_{T_S^J(x),\E}(f(x))\mid J \subset J(x)\},
\end{equation}
where $T_S^J(x) := \{v \in \R^n \, \vert \, \langle \nabla h_i(x), v \rangle \ge 0, i \in J\}$.
\end{proposition}

{\em Proof of Prop.~\ref{prop:kras}:}
As the set sequence $F(B(x,\delta))$ is monotonically increasing with $\delta$, i.e., $F(B(x,
\delta_1)) \subseteq F(B(x,\delta_2))$, if $\delta_1 \le \delta_2$,   \cite[Ex.~4.3(b), Prop.~4.30(b)]{rockafellar2} yield \eqref{eq:convKrasovskii}. % we can write $
%K_F(x) = {\rm con} \limsup_{y \to x} F(y),$
%where the lim sup on the right is interpreted in terms of convergence of sets. 
%For each $z \in \limsup_{y \to x} F(y)$, there exist  sequences $\{z_k\}_k$ and $\{y_k\}_k$ with $z_k \in \Pi_{S,\cE}(y_k,f(y_k))$ such that $z_k \to z$ and $y_k \to x$. Due to continuity of $f$, $f(y_k) \to f(x)$. Since $\|z_k - z\|$ gets arbitrarily small for large $k$, it follows that $z \in \lim_{k \to \infty} \Pi_{S,\cE}(y_k,f(x))$ \textcolor{red}{limsup??} and, hence, \eqref{eq:convKrasovskii} holds. 
%\textcolor{red}{maybe I am missing something obvious, but given that this is crucial step for (34) -- which is delicate, I want to be super careful: the result would imply something like: $\|\Pi_{S,\cE}(y_k,f(x))-\Pi_{S,\cE}(y_k,f(y_k)) \|$ converging to $0$. We only assume closedness of $S$, and later we apply this to sectors. certainly for non regular sets a small change in $f$ can result in different outcome, no? }

 ``$\supseteq$ \eqref{eq:krasovskiiFinite}:'' For $S$ as in \eqref{eq:defSFinite}, and a fixed $x \in S$,  
%\textcolor{red}{a subtlety here is that we use $T_S^{J(y)}\rightarrow T_S^{J(x)}$ when $y\rightarrow x$ [or at least the convergence of the projections on these sets], also in the last line of the proof this is used.} 
  suppose, w.l.o.g.,  $J(x) = \{1,\dots,m_x\}$. Note that we can choose vectors $v_1, \cdots, v_{m_x} \in \R^n$ such that the Jacobian of the mapping $ \R^{m_x} \ni \alpha:=(\alpha_1,\cdots,\alpha_{m_x}) \mapsto (h_1(x+ \alpha_1v_1 + \dots + \alpha_{m_x}v_{m_x}), \dots, h_{m_x}(x+ \alpha_1v_1 + \dots + \alpha_{m_x}v_{m_x})) \in \R^{m_x}$ is nonsingular at $\alpha= 0$ due to {\bf(CQ)}. Using the inverse function theorem, it follows that, for each $J \subset J(x)$ there exists a sequence $y_k \to x$ such that $J(y_k) = J$. Along this sequence, it holds that $\lim_{y_k \to x}P_{T_S(y_k),\cE}(f(y_k)) = \lim_{y_k \to x}P_{T_S^J(y_k),\cE}(f(y_k))= P_{T_S^J(x),\cE}(f(x))$, using \eqref{eq:Tsexplicit}. Now    ``$\supseteq$'' in   \eqref{eq:krasovskiiFinite} follows from \eqref{eq:convKrasovskii}.
 
 ``$\subseteq$:'' Consider a sequence $y_k \to x$  with  $z= \lim_{y_k \to x}P_{T_S(y_k),\cE}(f(y_k))$. Due to continuity of $h$,  we have $J(y_k)\subseteq J(x)$ for $k$ large enough. Therefore, as only a finite number $J(y_k)$'s are possible,  we can select a subsequence $y_{k_l} \to x$ with $z=\lim_{l \rightarrow \infty}P_{T_S(y_{k_l}),\cE}(f(y_{k_l}))$ and $J(y_{k_l}) = J \subset J(x)$ constant. Now  ``$\subseteq$'' in \eqref{eq:krasovskiiFinite} follows from the last step in  the ``$\supseteq$ \eqref{eq:krasovskiiFinite}''-proof, using characterisation \eqref{eq:convKrasovskii}. $\Box$

{\em Poof of Theorem~\ref{thm.kras=pds}:} The inclusion $\supseteq$ is obvious; so   consider the inclusion $\subseteq$. Thereto, let $v\in  K_F(x) \cap T_S(x)$. Based on Proposition~\ref{prop:kras}, $v\in K_F(x) $ can be written as 
    $ v= \sum_{J \subseteq J(x)} \lambda_J P_{T_S^J(x),\E}(f(x))$
    with $\lambda_J \geq 0$, for $J \subset J(x)$, and $\sum_{J \subset J(x)} \lambda_J=1$. Let us consider
    \begin{equation*}
    \begin{aligned}
        \|v-f(x)\| %&= \Big\| \sum_{J \subset J(x)} \lambda_J P_{T_S^J(x),\E}(f(x))-f(x)\Big\| \\
        & = \Big\|\sum_{J \subseteq J(x)} \lambda_J P_{T_S^J(x),\E}(f(x))-\sum_{J \subset J(x)} \lambda_J  f(x)\Big\| \\
        & \leq \sum_{J \subseteq J(x)} \lambda_J \big\| P_{T_S^J(x),\E}(f(x)) - f(x)\big\| \\
        & \leq \sum_{J \subseteq J(x)} \lambda_J \big\| P_{T_S^{J(x)}(x),\E}(f(x)) - f(x)\big\| \\
        & = \big\| P_{T_S(x),\E}(f(x)) - f(x) \big\|
         = \big\| \Pi_{S,{\E}}(x,f(x))-f(x) \big\|,
    \end{aligned}
    \end{equation*}
where we used in the second inequality that 
$\|P_{T_S^J(x),\E}(f(x)) - f(x)\| \leq \|P_{T_S^{J(x)}(x),\E}(f(x)) - f(x) \| $
as $T_S^{J(x)}(x) \subseteq T_S^J(x))$. 
Moreover, by definition of $P_{T_S^J(x),\E}$, it holds that $P_{T_S^J(x),\E}(f(x)) - f(x) \in \E$ for each $J\subset J(x)$ and thus,
\begin{equation*}
    \sum_{J \subseteq J(x)} \lambda_J (P_{T_S^J(x),\E}(f(x)) - f(x)) =  v - f(x) \in \cE %\\ \underbrace{\sum_{J \subseteq J(x)} \lambda_J P_{T_S^J(x),\E}(f(x))}_{=v} - f(x) \in \E
\end{equation*}
due to $\E$ being a linear subspace. Since $v\in T_S(x)$ and $v-f(x) \in \E$, and $\|\Pi_{S,\cE}(x,f(x)) - f(x) \|$ is the shortest distance along $\E$ between $T_S(x)$ and $f(x)$, i.e., 
%\[\Pi_{\cS,\E}(x,f(x)) = \argmin_{w \in T_S(x), w - f(x) \in \E} \|w-f(x)\|\]
it must hold that $\|v-f(x)\| = \|\Pi_{S,{\E}}(x,f(x))-f(x)\|$
 and thus $v$ must be the unique  closest point in $T_S(x)$ to $v$ along $\E$ and thus $v=\Pi_{S,{\E}}(x,f(x)) $, thereby proving the result. $\Box$

\begin{theorem}
    Assume $f$ is continuous and the set $S$ in \eqref{eq:defSFinite} satisfies {\bf (CQ)} and \eqref{eq:nonempty}. Then for every $x(t_0)=x_0 \in S$ there exist $T>0$ and an AC solution $x:[0,T] \rightarrow S$ to \eqref{pds2}. 
\end{theorem}

{\em Proof:}
     It follows from \cite[Lemma 5.16]{goebel2012hybrid} that $K_F$ is outer semicontinuous. Moreover, $K_F$ takes non-empty, closed and convex set-values, and for all $x\in S$ there is an open  neighborhood $U$ of $x$ such that for all $y\in U\cap S$ it holds that $K_F(y) \cap T_S(y) \neq \emptyset$ (as it contains $\Pi_{S,\E}(y, f(y))$). According to \cite[Lemma 5.26~(b)]{goebel2012hybrid}, see also \cite{aubin}, the corresponding viability conditions are satisfied implying that the Krasovskii regularization \eqref{eq:convKrasovskii} has a solution $x$. Hence, this solution $x$  satisfies $x(t)\in S$ for all $t\in [0,T]$, it holds that $\dot x(t) \in T_S(x(t))$, almost everywhere, see \cite[Lemma 5.26~(a)]{goebel2012hybrid}. Hence, it holds a.e. that $ \label{aux1} \dot x(t) \in K_S(x(t)) \cap T_S(x(t)).$  Invoking Theorem~\ref{thm.kras=pds} shows that $x$ is now a solution to \eqref{pds2}. $\Box$

\section{ePDS on Sectors} \label{sec:epds_sectors}

In the previous section we focussed on ePDS with partial projection, but regular, finitely generated sets. In this section, we abandon the regularity by allowing sectors as the constraint sets $\cS$, which is important for PBC. We start by showing that closed-loop PBC systems can be written in the form of ePDS, and prove that the corresponding projection operator is well-defined. In Section~\ref{sec:well} we provide results on existence and completeness of solutions.

\subsection{Closed-loop PBC systems are ePDS} \label{sec:motivation}

Consider the general nonlinear SISO plant given by 
\begin{subequations} \label{eq:plant5}
\begin{eqnarray}
\dot{x} & = & f_p(x, u) \\
e &= & G_px
\end{eqnarray} 
\end{subequations}
with state $x\in \ree^n$, control input $u\in \ree$ and output $e\in \ree$. Here, $f_p:\ree^{n}\times \ree \rightarrow \ree^n$ is a continuous function and $G_p\in \ree^{1 \times n}$ a row vector. This system is connected to a PBC for which the unprojected dynamics are given by 
\begin{subequations}\label{eq:controller5}
\begin{eqnarray}
\dot{z} &= & f_c(z,e)\\
u &= & z_1
\end{eqnarray} 
\end{subequations}
with state $z\in \ree^m$, controller output $u\in \ree$ and controller input $e\in \ree$. The map $f_c:\ree^{m}\times \ree \rightarrow \ree^m$ is assumed to be continuous. Note that we have for both plant and controller linear output equations, which in many cases  can be realized by suitable coordinate transformations.  The reason for this setting will become clear later. For ease of exposition, we took $u=z_1$.
%We assume that all the maps $f_p$, $g_p$ and $f_c$ are continuous. %Moreover, $g_p$ is assumed to be continuously differentiable.  {\bf check actual needs after proofs!}
The projection will take place only along controller states $z$-dynamics, as we cannot change the plant  dynamics as they adhere to physical laws, resulting in a {\em partial} projection operation with the goal to keep the output-input pair $(e,u)$ in a sector 
\begin{equation} \label{eq:sector5}
 S=   \left\{(e,u)\in \mathbb{R}^2 \mid (u-k_1e)(u-k_2e) \leq 0 \right\}, 
\end{equation}
where $k_1, k_2 \in \ree$ with $k_1<k_2$, as is motivated by HIGS and other PBCs \cite{DeeSha21,EijHee_CSL20a,ShiNi_CDC22,Deenen2017,ShiNi_CDC22}. 
Note that the set $S$ does not satisfy the constraint qualification {\bf (CQ)} in Section~\ref{sec:epds} as the gradient of the mapping, $\R^2 \ni (e,u)\mapsto (u-k_1e)(u-k_2e) \in \R$, vanishes at the origin, which also reflects the absence of Clarke regularity \cite{rockafellar2}.
We can also write $S$ as the union of two polyhedral cones 
\begin{subequations}
\begin{equation}
    S= K \cup -K \text{ with }
\end{equation}
\begin{equation} \label{eq:K} K=    \left\{(e,u)\in \ree^2 \mid u \geq k_1 e \text{ and } u \leq k_2 e \right\}. \end{equation}
\end{subequations}
To obtain the closed-loop system description, we introduce the state $\xi=(x,z) \in \ree^{n+m}$, the constraint set ${\cal S}\subseteq \ree^{n+m}$  as 
\begin{equation}\label{eq:defBigS}
 %{\cal S}=   \left\{\xi=(x,z)\in \mathbb{R}^{n+m} \mid \underbrace{(G_px, z_1)}_{=(e,u) =: H\xi}\in S \right\}, 
 {\cal S}=   \left\{\xi=(x,z)\in \mathbb{R}^{n+m} \mid (G_px, z_1)=(e,u) =: H\xi\in S \right\}, 
\end{equation}
and the projection subspace as 
\begin{equation}
    \E = \textup{Im}\underbrace{\begin{bmatrix} O_n \\ I_m \end{bmatrix}}_{=:E} \text{ and } H=\footnotesize{\begin{pmatrix} G_p & 0 \\ 0 & [ 1 \ 0 \ldots 0]\end{pmatrix}}, 
\end{equation} 
where $\textup{Im} E $ denotes the column space of the matrix $E$. 
%, where $e_{n+1}\in \ree^{m+n}$ denotes the $n+1$-th standard basis vector, i.e., the $n+1$-th element of $e_{n+1}$ is $1$, the others all $0$. 
%Note that we can also write $\cS = \cK_1 \cup \cK_2 $, 
%where $\cK_i$, $i=1,2$, can be put in the form \eqref{eq:defSFinite}, particularly, 
%$\cK_1=\{\xi \mid h_i(\xi) \geq 0, \ i=1,2\}$ and 
%$\cK_2=\{\xi \mid -h_i(\xi) \geq 0, \ i=1,2\}$ with
%\begin{subequations}
%\begin{eqnarray}
%    h_1(\xi) &= & u-k_1e = z_1 - k_1 g_p(x_p),\\
%    h_2(\xi) &= & -u+k_2e = -z_1 + k_2 g_p(x_p).
%\end{eqnarray}    
%\end{subequations}
The closed-loop dynamics can now be written as  the ePDS
\begin{equation} \label{eq:cloop}
    \dot \xi = \Pi_{{\cal S},\E}(\xi,f(\xi)) := F(\xi)
\end{equation}
with $f$ denoting the unprojected closed-loop vector field
\begin{equation}
    f(\xi) = (f_p(x, z_1),f_c(z,G_p x)), 
\end{equation} which we sometimes write, with some abuse of notation, as $(f_p(\xi),f_c(\xi))$, and we split $f_c(\xi)$ as $(f_{c,1}(\xi),f_{c,2}(\xi))$ with $f_{c,1}(\xi) \in \ree$ and thus $\dot u =  \dot z_1= f_{c,1}(\xi)$. Hence, $f_{c,2}(\xi) \in \ree^{m-1}$.

\subsection{Well-defined partial projection operators on sectors}

Consider the ePDS \eqref{eq:cloop} with sector constraints as in the previous subsection. As a first step towards establishing existence of solutions, we  show that the introduced partial projection $ \Pi_{S,E}(x,v)$ provides a unique outcome for each $x\in S$ and each $v\in \ree^n$, for which the next lemma is useful. %To do so, we introduce the following result.  

\begin{lemma} \label{lem:tc}
    Let $D\subseteq \ree^{n_d}$ be a closed set and $H\in\ree^{n_c \times n_d}$ with $H$ full row rank. Define $C\subseteq \ree^{n_c}$ as $C:=\{c \in \ree^{n_c}\mid Hc \in D\}$. Let $x\in C$. Then $T_C(x)= \{v \in \ree^{n_c} \mid  Hv \in T_D(Hx)\}$.
\end{lemma}
%\begin{proof}
%    to be added
%\end{proof}

{\em Proof:} ``$\subseteq$'' follows from \cite[Thm.~6.31]{rockafellar2}. %$v\in T_C(x)$, then there are sequences $\{c_i\}_{i\in\N} \subseteq C$ and $\{t_i\}_{i\in\N} \subseteq \ree_{> 0}$ with  $x=\lim_{i \rightarrow \infty}c_i$,  $\lim_{i \rightarrow \infty} t_i=0$ and $v=\lim_{i \rightarrow \infty}\frac{c_i-x}{t_i}$. Clearly, $Hv=\lim_{i \rightarrow \infty}\frac{Hc_i-Hx}{t_i}$, and $\{Hc_i\}_{i\in\N} \subseteq D$ with $Hx=\lim_{i \rightarrow \infty}Hc_i$. Hence, $Hv \in T_D(Hx)$.
``$\supseteq$:'' Consider $v$ with $Hv \in T_D(Hx)$. Hence, there are sequences $\{d_i\}_{i\in\N} \subseteq D$ and $\{t_i\}_{i\in\N} \subseteq \ree_{> 0}$ with  $Hx=\lim_{i \rightarrow \infty}d_i$,  $\lim_{i \rightarrow \infty} t_i=0$ and $Hv=\lim_{i \rightarrow \infty}\frac{d_i-Hx}{t_i}$. We will use the decomposition $\R^{n_c} = \ker H \oplus \im H^\top$, where $\ker H = \{ w\mid Hw=0\}$, and the projection $H^\top(HH^\top)^{-1}H$ on $\im H^\top$ along $\ker H$. We write $x$ uniquely as $x=x_a+x_b$ with $x_a\in \ker H$ and $x_b=H^\top(HH^\top)^{-1}Hx \in \im H^\top$ and, similarly $v=v_a+v_b$ with $v_a \in \ker H$ and $v_b=H^\top(HH^\top)^{-1}Hv$. Take $c_i:=x_a+ H^\top(HH^\top)^{-1}d_i + t_i v_a$. Note that $Hc_i=d_i \in D$ and thus $c_i\in C$. Moreover, $c_i$ converges to $x_a+ H^\top(HH^\top)^{-1}Hx=x_a+x_b=x$, when $i\rightarrow \infty$. Finally, note  $\frac{c_i-x}{t_i} = \frac{H^\top (HH^\top)^{-1}d_i -x_b}{t_i} + v_a$, which is equal to $H^\top (HH^\top)^{-1}(\frac{d_i -Hx}{t_i})+ v_a \longrightarrow H^\top (HH^\top)^{-1}Hv+ v_a =v_b + v_a=v$. Hence, $v \in T_c(x)$. $\Box$

This result will be instrumental below and explains the choice for the linear output equations for  plant \eqref{eq:plant5} and controller  \eqref{eq:controller5}. In case the output equations would be nonlinear, an extension of Lemma~\ref{lem:tc} would be needed, which due to the absence of Clarke regularity of the sector \eqref{eq:sector5} in the origin,  is not straightforward, see, e.g.,   \cite[Thm.~6.31]{rockafellar2}. Extending this lemma to nonlinear maps is interesting  future work. 

Exploiting Lemma~\ref{lem:tc} gives 
\begin{subequations} \label{eq:to2d}
\begin{align}
 & \Pi_{{\cal S},\E}(\xi,f(\xi)) %=  \\ 
   =  \argmin_{v \in T_{\cal S}(\xi), f(\xi)-v \in \E} \|f(\xi)-v\|  \\
   &= \argmin_{v=(v_x,v_z) \in T_{\cal S}(\xi), v_x=f_p(\xi)} \|f_c(\xi)-v_z\|  \\ 
    &=  \argmin_{v=(f_p(\xi),v_z) \in T_{\cal S}(\xi)} \|f_c(\xi)-v_z\|  \\ 
     & \stackrel{Lem.~\ref{lem:tc}}{=}  \argmin_{v=(f_p(\xi),v_z), (G_p f_p(\xi), v_{z,1})  \in T_{S}(G_p x,z_1)} \|f_c(\xi)-v_z\|  \\ 
      &= \argmin_{v=(f_p(\xi),v_{z,1},f_{c,2}(\xi)), (G_p f_p(\xi), v_{z,1})   \in T_{S}(G_p x,z_1)} |f_{c,1}(\xi)-v_{z,1}|.
\end{align} 
\end{subequations}
Hence, this shows that 
\begin{subequations} \label{eq:structurePi}
\begin{equation} 
\Pi_{{\cal S},\E}(\xi,f(\xi)) = (f_p(\xi),v_{z,1}^*,f_{c,2}(\xi)) \text{ with }\end{equation}
\begin{equation} \label{eq:vz1*}
    v_{z,1}^* = \argmin_{v_{z,1}, (G_p f_p(\xi), v_{z,1})   \in T_{S}(G_p x,z_1)} |f_{c,1}(\xi)-v_{z,1}|.
\end{equation}
\end{subequations}
Note that \eqref{eq:to2d} reveals that the partial projection  $\Pi_{{\cal S},\E}(\xi,f(\xi)) $ only alters the $u=z_1$-dynamics and the rest remains unchanged, including the $e$-dynamics. In fact, the projection $\Pi_{{\cal S},\E}$ can be perceived to take place in a 2-dimensional subspace. Indeed, the optimization in \eqref{eq:vz1*} can be retrieved as (part of) a partial projection in a 2-dimensional $(e,u)$-space. To make this more concrete, the unprojected $(e,u)$-dynamics can be written as 
\begin{equation}
    (\dot e,\dot u) = (G_p f_p(\xi),f_{c,1}(\xi))=: f_{eu}(\xi)
\end{equation} and the corresponding projected dynamics
\begin{equation} \label{eq:fromndto2d}
  (\dot e,\dot u) = \Pi_{S,\cE'}((e,u),  f_{eu}(\xi)) = (G_p f_p(\xi),v_{z,1}^*)
\end{equation}
 with $\cE'=\im E'$ where $E'= {\footnotesize{\begin{pmatrix}0\\1 \end{pmatrix}}}$. 
Interestingly, $\Pi_{S,\cE'}(s,w)$ for  $s=(e,u)\in \ree^2$ and $w=f_{eu}(\xi)\in \ree^2$ can be written as
\begin{equation} \label{piSE2}
\Pi_{S,\cE'}(s,w) = w + E' \eta^*(s,w),
\end{equation}
where
\begin{equation} \label{piSE2a}
 \eta^*(s,w) = \argmin_{\eta \in \Lambda(s,w)} \underbrace{\|E'\eta\|}_{=|\eta|} \text{ and }
\end{equation}
\begin{equation} \label{piSE2b}
\Lambda(s,w) = \{ \eta \in \ree \mid  w+E'\eta \in T_S(s)\}
\end{equation}
 
Using the above  and observing that 
\begin{equation} \label{eq:SE}
   E' \text{ full column rank }, \cE' \cap S = \{0\} \text{ and } S + \cE' = \ree^{2},
\end{equation}
we can obtain a lemma that will lead to well-definedness of $\Pi_{\cS,\cE}$ and $\Pi_{S,\cE'}$.

\begin{lemma} \label{lem:wp} Consider $S$ as in \eqref{eq:sector5} for $k_2 > k_1$ and let $\cE'=\im E'$ satisfying \eqref{eq:SE}. Then, it holds 
 for each $s\in S$ and each $w\in \ree^2$  that $\Lambda(s,w)$ is a non-empty  closed polyhedral set. 
\end{lemma}

{\em Proof:}
Let $s\in S$ and  $w\in \ree^2$ be given and notice that 
\begin{equation} \label{eq:TSx}
T_S(s) = \begin{cases} 
T_K(s), & \text{ when } s\in K\setminus -K,\\
K\cup - K, & \text{ when } s\in K \cap -K,\\
-T_K(s), & \text{ when } s\in -K\setminus K.
\end{cases}
\end{equation}
It follows from  $S+ \cE' = \ree^2$ that for all $s\in 
S$ we have $T_S(s) + \cE' \neq \emptyset$ and thus $\Lambda(s,w) $ is non-empty. Clearly, when $s\not\in K \cap -K$ it follows that $T_S(s)$, as given in \eqref{eq:TSx}, is a closed polyhedral cone (as $K$ is polyhedral cone, see \eqref{eq:K}) and then  $\Lambda(s,w) $ is a closed polyhedral set. 
 So, let us focus on $s\not\in K \cap -K$, where $T_S(s)=K\cup -K$ and thus $\Lambda(s,w) = \{ \eta \in \ree \mid  w+E'\eta \in K\cup -K\}$. 

\underline{Claim}: $w+ E'\eta \in K$ and $w+ E'\bar \eta \in -K$ imply that $\eta = \bar \eta$. 
To prove the claim, note that due to $K$ being a convex cone and $-w- E'\bar \eta \in K$, we get that 
\[ E'(\eta - \bar \eta) = (w+ E'\eta) -w- E'\bar \eta \in K\]
Since $\cE'\cap S = \{0\}$ and $E'$ has full column rank, this shows that $\eta=\bar \eta$ and the result follows. 

Using the Claim, it follows that if there is $\eta$ with $w+ E'\eta \in K$ then $\Lambda(s,w)  = \{ \eta \in \ree\mid  w+E'\eta \in K\}$ (as any $\eta$ in $\{ \eta \in \ree\mid  w+E'\eta \in -K\}$ would also be contained in $ \{ \eta \in \ree\mid  w+E'\eta \in K\}$, and similarly if there is $\bar \eta$ with $w+ E'\bar \eta \in -K$ then $\Lambda(s,w)  = \{ \eta \in \ree\mid  w+E'\eta \in -K\}$. As the sets $ \{ \eta \in \ree\mid  w+E'\eta \in  K\}$  and $ \{ \eta \in \ree\mid  w+E'\eta \in - K\}$   are both closed polyhedral sets, so is $\Lambda(s,w) $. $\Box$

Due to the constraint set of \eqref{piSE2a} being a closed polyhedral set, and the square of the cost function of \eqref{piSE2} is $\eta^\top (E')^\top E' \eta = \eta^2$ is a quadratic positive-definite function, a unique minimizer exists, showing that \eqref{eq:fromndto2d} is well-defined and so is \eqref{eq:structurePi}, i.e., $\Pi_{{\cal S},\E}(\xi,v)$ for all $\xi\in \cS$ and $v\in \ree^{m+n}$.

Interestingly, from the proof of Lemma \ref{lem:wp} and \eqref{eq:TSx}, it follows that 
for all $s\in S$ and $w\in \ree^{2}$
\begin{equation} \label{eq:PSinPK'}
\Pi_{{S},\E'}(s,w) \in  \left\{ \Pi_{{ K},\E'}(s,w),\Pi_{-{K},\E'}(s,w) \right\}
\end{equation}
and, using Lemma~\ref{lem:tc}, we obtain for all $\xi\in \cS$, $v\in \ree^{m+n}$
\begin{equation} \label{eq:PSinPK}
\Pi_{{\cal S},\E}(\xi,v) \in  \left\{ \Pi_{{\cal K},\E}(\xi,v),\Pi_{-{\cal K},\E}(\xi,v) \right\},
\end{equation}
with the convention that some of the projections in the set in the right-hand sides of \eqref{eq:PSinPK} and \eqref{eq:PSinPK'} may be empty, e.g., $\Pi_{{\cal K},\E}(\xi,v)=\emptyset$, if $(v+ \E) \cap T_K(\xi) = \emptyset$ (which is, amongst others, the case when $\xi\not\in \cK$).

\section{Well-posedness of ePDS on Sectors} \label{sec:well}

We focus now on existence and completeness of solutions.

\begin{theorem} \label{thm:pbc} Consider \eqref{eq:cloop} with set $\cS$ defined via \eqref{eq:defBigS} and \eqref{eq:sector5}. 
    For each initial state $\xi_0 \in \cS$, there exists a Carath\'eodory solution locally, i.e., there is a $T>0$ such that $\xi:[0,T]\rightarrow \ree^{n+m}$ with $\xi(0)=\xi_0$ is a solution to \eqref{eq:cloop}.
\end{theorem}

%In the proof of Theorem~\ref{thm:pbc}, we use the Krasovskii regularisation of \eqref{eq:cloop} given by 
%$\dot x \in \cap_{\delta >0} \overline{\textup{con}}F(B(x,\delta)) =: K_F(x)$, as before. {\bf I think we can generally write under continuity of $f$:}
%\begin{equation} \label{eq:kras}
%   {\color{red}
%K_F(\xi) = \overline{\rm con} \limsup_{\bar \xi \to \xi} P_{T_\cS(\bar \xi),\cE} (f(\xi))}
%\end{equation}
%{\bf this one is crucial/ makes life easier as due to the structure of $\cS$ and $S$ we have finitely many here -- most insightful for our proofs would be to have here 2d version linked to \eqref{eq:fromndto2d}}

%\begin{proposition} \label{prop:kras_cloop}
%The Krasovskii regularization of $F$ at $\xi\in {\cal S}$ is equal to  {\bf to be completed:}
%\begin{equation}
 %   K_F(x) = \textup{con}\{P_{T_J,\E}(f(x))\mid J \subset J(x)\}
%\end{equation}
%\end{proposition}
%\begin{proof}
 %  Using section 2/3, we can exploit now \eqref{eq:PSinPK} to %immediately get this counterpart!!! note that $f$ continuous. 
%\end{proof}

\begin{figure}[b!]
\centering
\begin{tikzpicture}[thick, xscale=0.75,yscale=0.85]

\coordinate (origin) at (0,0);
\coordinate (ltop) at (-3cm,0);
\coordinate (lbot) at (-3cm,-2cm);
\coordinate (rbot) at (3cm,0);
\coordinate (rtop) at (3cm,2cm);

\fill[black!20] (ltop) -- (origin) -- (lbot) -- cycle;
\draw[very thick, black] (ltop) -- (origin) -- (lbot);

\fill[black!20] (rbot) -- (origin) -- (rtop) -- cycle;
\draw[very thick, black] (rbot) -- (origin) -- (rtop);

\draw[Darkred, very thick, ->] (origin) -- ++(340:2.15cm) node[anchor = west] {$f_0 = f(\xi)$};
\draw[Darkred, very thick, ->] (origin) -- ++(0:2cm) node[anchor = south west] {$f_1 = \Pi_{S,\cE}(f(\xi))$};
\draw[Darkred, very thick, ->] (origin) -- ++(34:2.35cm) node[anchor = south east] {$f_2 $};
\draw[UIUCorange, very thick, dashed] (2,-2) -- (2,2.2) node[anchor = north west]{$\cE$};
\draw[black, thick] (-1,-0.5) node[anchor=east]{$\cS$};

\end{tikzpicture}
\caption{Depiction of a counterexample showing that $K_F(\xi) \cap T_S(\xi) \neq \Pi_{S,\cE}(\xi,f(\xi))$, for $\xi \in \cK \cap -K$. In particular, $K_F(\xi) = {\rm co}\{f_0,f_2\}$ and $K_F(\xi) \cap T_S(\xi) = {\rm co}\{f_1,f_2\} \neq \{f_1\}$.\label{fig:counterexample}}
\label{fig:counterEx}
\end{figure}
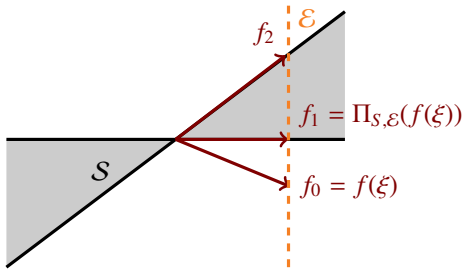

The proof is based, as in Section~\ref{sec:epds}, on showing existence of Krasovskii solutions, and then proving that the Krasovskii solutions are also (Carath\'eodory) solutions to  \eqref{eq:cloop}. Proving Theorem~\ref{thm:pbc} for sectors is more complicated compared to the counterpart in Section~\ref{sec:epds} due to absence of (Clarke) regularity \cite{rockafellar2}. For this reason, the relation \eqref{eq:viab=PDS} does not hold everywhere for sectors as illustrated in Fig.~\ref{fig:counterexample}. However, we show that Krasovskii solutions will only visit states for which \eqref{eq:viab=PDS} is violated on a set of times with measure zero. 

%\begin{figure}[!h]
%\centering
%\includegraphics[width=7cm]{KintTnotPi}
%\caption{Due to the absence of convexity/prox regularity, we do not have that $K_F(\xi) \cap T_\cS(\xi) = \{\Pi_{\cS,\E}(\xi, f(\xi))\}$. This is due to the tangent cone lacking a lower semicontinuity property (in the set $\cK \cap -\cK$.)}
%\end{figure}

{\em Proof Thm.~\ref{thm:pbc}:} First note that the Krasovskii regularisation $K_F$ is outer semicontinuous and takes non-empty convex closed set-values. Moreover, observe that $K_F(\xi) \cap T_\cS(\xi)\neq \emptyset$ for all $\xi \in S$ as $\Pi_{{\cal S},\E}(\xi,f(\xi))$ is contained in the intersection. Hence, the local viability assumption, see \cite[Lemma 5.26]{goebel2012hybrid} and \cite{aubin}, is satisfied and we can establish for any initial state the existence of an AC solution $\xi:[0,T] \to \cS$ to $\dot \xi\in K_F(\xi)$ for some $T>0$. Clearly, due to the necessity of the viability condition in \cite[Lemma 5.26]{goebel2012hybrid}, we obtain that the solution satisfies, for almost all times $t\in [0,T]$, that 
$ \dot \xi \in K_F(\xi) \cap T_{\cal S}(\xi).$
Interestingly, for all times $t\in [0,T]$ for which $\xi(t) \not \in K\cap -K$, we can use the results in Theorem~\ref{thm.kras=pds}, as 
$
    \Pi_{{\cal S},\E}(\xi,f(\xi)) = \Pi_{{\cal K},\E}(\xi,f(\xi))
$
in a neighbourhood of $\xi \in \cK \setminus -\cK$ (and similarly for $\xi \in -\cK \setminus \cK$). 
Hence,  for almost all $t$ with $\xi(t) \not \in K\cap -K$ (and due to the continuity of $\xi$ and closedness of $K\cap -K$, there is $\epsilon_t$ such that for $\tau \in [t,t+\epsilon_t)$, we have $\xi(\tau) \not \in K\cap -K$) 
\begin{equation} \label{aux4c} K_F(\xi(t)) \cap T_S(\xi(t)) = \Pi_{S,\E}(\xi(t), f(\xi(t))).\end{equation} 
Hence, we only have to consider the times where $\xi(t)  \in K\cap -K$, i.e. $(e(t),u(t))=(0,0)$.
We will consider two cases: (i) $\dot e(t) \neq 0$, and
(ii) $\dot e(t) = 0$, and only  times $t$ where $\xi$ is differentiable. The latter can be done as the set of non-differentiability has measure zero.

Case (i): w.l.o.g.,  assume $\dot e(t) > 0$ %(when $\dot e(t) < 0$ the arguments are similar). 
 Note that $\dot e = G_p f_p(\xi)$, which is a continuous function of time along solution $\xi$ (as $f_p$ and $\xi$ continuous). Hence,  there are $\epsilon_t$ and $\eta>0$ such that $\dot e(\tau) \geq \eta $ for $\tau \in [t,t+\varepsilon_t]$. Thus,   $e(\tau)\geq \eta (\tau- t)>0$ for $\tau \in (t,t+\epsilon_t]$. Hence, the time $t$ where $e$ is zero is an isolated point  in this case, and for $ \tau \in (t,t+\epsilon_t]$ $\xi(t) \in \cK \setminus -\cK$, for which  \eqref{aux4c} holds (with $t=\tau$).  In fact, the set $\{t \in [0,T]\mid e(t)=u(t)=0 \text{ and }\dot e(t) \neq 0 \}$ is countable as its points are right-isolated, and, hence, in each interval $(t,t+\varepsilon_t]$ there lies a rational number not in this set, and the rational numbers in $[0,T]$ are countable. Therefore, this set is of measure zero. 

Case (ii) $\dot e(t) = G_p f_p(\xi(t))= 0$. 
%{\color{red}Note that $(\dot e(t),\dot u(t)) \in T_{\bar S}(0,0)=S$. %We exploit now that, due to the structure of the projection operator $\Pi_{\cS,\cE}$ as revealed in \eqref{eq:structurePi} and \eqref{eq:fromndto2d},  and continuity of $f_p$ and $f_c$ that the Krasovskii regularization can be written for all $\xi \in \cS$ as 
%\begin{equation} \label{eq:Kform}
 %   K_F(\xi) = \{f_p(\xi)\} \times Z_1(\xi) \times \{f_{c,2}(\xi)\} \end{equation}
%for some osc mapping $Z_1: \cS \rightarrow \ree$. In particular, for case (ii), where  $T_{S}(0,0)=S$, and 
%
%As $\dot e(t)= G_p f_p(\xi(t))=0$, we obtain from the ``2-dimensional perspective'' \eqref{eq:fromndto2d} and \eqref{eq:convKrasovskii}  combined with \eqref{eq:PSinPK'}-\eqref{eq:PSinPK}, and the characterisation in Prop.~\ref{prop:kras} (note that in the $(e,u)$-space the direction of $f$ is $(0,f_{c,1}(\xi))$, compare this to the sector $S$), that
 %\begin{equation} \label{eq:x5}
%    K_F(\xi) = \textup{con}\{f(\xi),(f_p(\xi),0,f_{c,2}(\xi))\}. 
%\end{equation} In \eqref{eq:x5} we used that $\Pi_{\cS,\cE}(\bar \xi, f(\xi)) =(f_p(\xi),0,f_{c,2}(\xi)) $ when $T_S(\bar \xi) \neq \ree^{m+n}$ for $\bar \xi \in \cS$. {\bf Maurice, shall we replace red with blue?}} 
 Now $f(\xi)$ takes the form $(0,f_{c,1}(\xi))$ in $(e,u)$-space. % and  \eqref{eq:fromndto2d} yields that  for $\bar \xi$ in a neighborhood of $\xi$ that either  $\Pi_{\cS,\cE}(\bar \xi, f(\xi)) =(f_p(\xi),0,f_{c,2}(\xi))$ or $\Pi_{\cS,\cE}(\bar \xi, f(\xi)) = f(\xi)$. 
 From \eqref{eq:convKrasovskii}, using \eqref{eq:structurePi}, \eqref{eq:fromndto2d}, we get
\begin{equation} \label{eq:x5}
    K_F(\xi) = \textup{con}\{f(\xi),(f_p(\xi),0,f_{c,2}(\xi))\}
\end{equation}
by considering all possible tangent cones $T_S(s) $ for $s\in S$ in a neighborhood of $(0,0)$. 
We either have $f(\xi) = (f_p(\xi),0,f_{c,2}(\xi))$ in which case  $K_F(\xi)$ is a singleton and thus must be equal to $\{\Pi_{{\cal S},\E}(\xi,f(\xi))\}$ (as this one is guaranteed to lie in $K_F(\xi)$), or $f(\xi) \neq (f_p(\xi),0,f_{c,2}(\xi))$ (so, $f_{c,1}(\xi) \neq 0$). In the latter case, $f(\xi) \not\in T_{\cal S}(\xi)$. Given the structure of $T_{\cal S}(\xi)={\cal S}$ (as $H\xi=(0,0)$) and $(1-\alpha)(f_p(\xi),0,f_{c,2}(\xi))+\alpha f(\xi)\not\in T_{\cal S}(\xi)$ for $\alpha \in (0,1)$ -- use here that $T_{\cS}(\xi) = \{v \in \ree^{m+n} \mid Hv \in T_S(H\xi)=S\}$ due to Lemma~\ref{lem:tc} -- it follows that  
\begin{equation}
    K_F(\xi) \cap T_{\cal S}(\xi)= \{(f_p(\xi),0,f_{c,2}(\xi))\}= \{\Pi_{{\cal S},\E}(\xi,f(\xi))\}
\end{equation}
Hence, we have $ \dot \xi(t)= \Pi_{S,\E}(\xi(t), f(\xi(t)))$ (if $\dot \xi(t)$ exists) in this case. Summarizing all cases, we obtain  a.e. $ \dot \xi(t)= \Pi_{S,\E}(\xi(t), f(\xi(t)))$. Hence,  $\xi:[0,T]\rightarrow \cS$ is  a solution.  $\Box$

Theorem~\ref{thm:pbc} shows {\em local} existence of Carath\'eodory solutions for a given initial condition. Below we extend this result to the existence of {\em global} solutions under suitable boundedness conditions of the unprojected vector field $f$.

\begin{corollary}\label{cor:completeness} Consider \eqref{eq:cloop} with sectorset $S$ as in \eqref{eq:sector5}. Moreover, we assume that there is $M>0$ such that  $\|f(\xi)\| \leq M (1+\|\xi\|)$ for all $\xi \in \cS$. For each initial state $\xi_0 \in \cS$, there exists a  Carath\'eodory solution $\xi:[0,\infty) \rightarrow \ree^{n+m}$ with $\xi(0)=\xi_0$ to \eqref{eq:cloop} on $[0,\infty)$.
\end{corollary}

{\em Proof:}
The proof starts by showing that the bound $\|f(\xi)\| \leq M (1+\|\xi\|)$ leads to a similar bound on $\|\Pi_{{\cal S},\E}(\xi,f(\xi))\|$ in \eqref{eq:cloop}.
To show this, recall \eqref{eq:structurePi} and \eqref{eq:vz1*}. Clearly, the bounds carry over to the $f_p(\xi)$- and $f_{c,2}(\xi)$-parts in \eqref{eq:structurePi}, so we only have to show that a similar bound also applies to $|v_{z,1}^*|$. Thereto, realize that \eqref{eq:PSinPK'} implies that $v_{z,1}^*$ satisfies
\begin{equation} \label{eq:vz1*'}
    v_{z,1}^* = \argmin_{v_{z,1}, (\dot e, v_{z,1})   \in T_{K}(e,z_1)} |f_{c,1}(\xi)-v_{z,1}|,
\end{equation}
where we replaced $G_p f_p(\xi)$ by $\dot e$,  $G_p x$ by $e$ and  $f_{c,1}(\xi)$ by $f_{c,1}$ for shortness, (or $v_{z,1}^*$ is given by  the same expression with $K$ replaced by $-K$). Using the form of $K$ in \eqref{eq:K} and consider all variations for $T_{K}(e,z_1)$ using the explicit expression in \eqref{eq:Tsexplicit}, we get that $(\dot e, v_{z,1})\in T_{K}(e,z_1)$ is equivalent to either $v_{z,1}\in \ree$, $k_1 \dot e\leq v_{z,1}$, $v_{z,1}\leq k_2 \dot e$, or $k_1 \dot e\leq v_{z,1}\leq k_2\dot e$, depending on the active constraint set in \eqref{eq:activeconstraints} for $K$. In fact, we obtain a piecewise linear solution for   $v_{z,1}^*$ in terms of $\dot e$ and $ f_{c,1}$, for each of the four options for $T_{K}(e,z_1)$, taking values  $v_{z,1}^*\in \{
      f_{c,1}, k_1 \dot e,k_2 \dot e\}$. Clearly, this shows that $|v_{z,1}^*| \leq \max(|f_{c,1}(\xi)|, |k_1| |G_p f_p(\xi)|, |k_2| |G_p f_p(\xi)|) \leq c \|f(\xi)\|$ for some $c>0$, and combining this with the bound on $f$, yields the existence of $M'>0$ such that also for all $\xi \in \cS$ we have 
     \[ \|\Pi_{{\cal S},\E}(\xi,f(\xi))\| \leq M' (1+ \|\xi\|).\]
   Using this bound and local existence of solutions per Theorem~\ref{thm:pbc}, we can now proceed similarly as in the proof of \cite[Thm.~4.2]{DeeSha21} to show by contradiction that a  maximal solution (i.e., a solution defined on the largest interval of the form $[0,T]$ possible) must be complete (i.e., $T=\infty$).  $\Box$

Consider now system \eqref{eq:cloop} with the inclusion of external time-varying functions in the plant, i.e., $\dot x = f_p(x,u,w)$ and $w$ a piecewise continuous function, i.e., $w\in PC$, meaning that there is 
$\{t_k\}_{k\in \N}\subset [0,\infty)$ with $t_0=0$, $t_{k+1}>t_k $ for all $k\in \N$, $\lim_{k\rightarrow \infty} t_k = \infty$,  $w$ is continuous for all $t\not\in \{t_k\}_{k\in \N}$, and  $\lim_{t \downarrow t_k} w(t) = w(t_k)$, $k\in\N$. Using the same controller and the modelling as in Section~\ref{sec:motivation}, we obtain the closed loop
\begin{subequations}\label{eq:cloop_input_all}
\begin{equation} \label{eq:cloop_input}
    \dot \xi = \Pi_{{\cal S},\E}(\xi,f(\xi,w(t))) \text{ with }
\end{equation}
\begin{equation}
    f(\xi,w) = (f_p(x, z_1,w),f_c(z,G_p x)), 
\end{equation}
\end{subequations}

\begin{corollary}
Consider \eqref{eq:cloop_input_all} with  $\cS$ as in \eqref{eq:defBigS} and $f$ continuous. Then for each $w\in PC$ and $\xi_0 \in \cS$, there is a $T>0$ such that a Carath\'eodory solution exists on $[0,T]$ to \eqref{eq:cloop_input_all} with $\xi(0)=\xi_0$ and input $w$.
Moreover, if for an $M>0$ 
\begin{equation} \label{eq:bound_cor}\|f(\xi,w)\| \leq M (1+\|(\xi,w)\|) \text{  for all } \xi \in \cS,
\end{equation} then for each initial state $\xi_0 \in S$ and bounded $w\in PC$ there exists a  Carath\'eodory solution $\xi:[0,\infty) \rightarrow \ree^{n+m}$ with $\xi(0)=\xi_0$  and input $w$ to \eqref{eq:cloop_input_all} on $[0,\infty)$.
\end{corollary}

{\em Proof:}
The idea of the proof is to embed $t$ as a state in $\chi = (\xi,t)$, see also \cite[p.~191]{aubin}, leading to the ePDS model \begin{equation} \label{eq:model_t}
    \dot \chi = \Pi_{\tilde{{\cal S}},\tilde{\cE}}(\chi,\tilde{ f}(\chi) )
\end{equation}
with $\tilde{{\cal S}} = \{ (\xi,t) \mid t \geq 0 \text{ and } \xi \in \cS\}$, $\tilde{ f}(\chi) = (f(\xi,w(t)), 1)$ and $\tilde{\cE}= \cE \times \{0\}$. Obviously, a solution $\xi$ to \eqref{eq:cloop_input} given input $w$ and initial state $\xi(0)=\xi_0$ leads to a solution $t\mapsto (\xi(t),t)$ to \eqref{eq:model_t} with $\chi(0)=(\xi_0,0)$ (without input) and vice versa. On $[0,t_1]$ $w$ is continuous, implying that $\tilde{ f}$ is a  continuous function of $\chi$ on $[0,t_1]$. Applying Theorem~\ref{thm:pbc} to \eqref{eq:model_t} proves now the local existence of solutions, and applying Corollary~\ref{cor:completeness} using  the bound \eqref{eq:bound_cor} guarantees that the solution is defined on the full interval $[0,t_1]$. Exploiting the bound \eqref{eq:bound_cor} again, following similar steps as in the proof of \cite[Thm.~4.2]{DeeSha21}, we obtain that the left limit $\lim_{t \uparrow t_1} 
\chi(t)$ exists; let us call this left limit $\chi(t_1)$, which lies in $\tilde{{\cal S}}$. Now we can repeat the arguments for the time window $[t_1,t_2)$ and, in fact, for each window $[t_{k},t_{k+1})$ leading to a solution on  $[0,t_k)$ for each $k\in\N$ and as $t_k\rightarrow \infty$ when $k\rightarrow \infty$, this leads to a solution on $[0,\infty)$.  $\Box$

\section{Conclusions}
We established existence and completeness results for solutions to extended PDS and closed-loop PBC systems (with and without inputs). This required careful analysis due to  partial projection operation and irregular constraint sets (sectors), which are important in PBC systems \cite{DeeSha21,EijHee_CSL20a,ShiNi_CDC22,Deenen2017,ShiNi_CDC22}, but obstructed the application of existing results. The results  provide cornerstones for further analysis  of PBC systems. 

\bibliographystyle{plain}
\bibliography{bib_nahs_v3}

\begin{thebibliography}{10}

\bibitem{aubin}
J.~P. Aubin and A.~Cellina.
\newblock {\em Differential Inclusions}.
\newblock Springer-Verlag Berlin Heidelberg, 1984.

\bibitem{brogliatoSCL}
B.~Brogliato, A.~Daniilidis, C.~Lemar\'echal, and V.~Acary.
\newblock On the equivalence between complementarity systems, projected systems
  and differential inclusions.
\newblock {\em Syst. Control Lett.}, 55:45--51, 2006.

\bibitem{Deenen2017}
D.A. Deenen, M.F. Heertjes, W.P.M.H. Heemels, and H.~Nijmeijer.
\newblock {Hybrid integrator design for enhanced tracking in motion control}.
\newblock In {\em 2017 American Control Conference}, pages 2863--2868. IEEE,
  2017.

\bibitem{DeeSha21}
D.A. Deenen, B.~Sharif, S.J.A.M. van~den Eijnden, H.~Nijmeijer, W.P.M.H.
  Heemels, and M.F. Heertjes.
\newblock {Projection-Based Integrators for Improved Motion Control:
  Formalization, Well-posedness and Stability of HIGS}.
\newblock {\em Automatica}, 133:109830, 2021.

\bibitem{dupuis1993dynamical}
P.~Dupuis and A.~Nagurney.
\newblock Dynamical systems and variational inequalities.
\newblock {\em Annals of Operations Research}, 44(1):7--42, 1993.

\bibitem{goebel2012hybrid}
R.~Goebel, R.G. Sanfelice, and A.R. Teel.
\newblock {\em Hybrid Dynamical Systems: modeling, stability, and robustness}.
\newblock Princeton University Press, 2012.

\bibitem{hauswirth2021projected}
A.~Hauswirth, S.~Bolognani, and F.~Dörfler.
\newblock Projected dynamical systems on irregular, non-euclidean domains for
  nonlinear optimization.
\newblock {\em SIAM Journal on Control and Optimization}, 59(1):635--668, 2021.

\bibitem{nagurney2012projected}
A.~Nagurney and D.~Zhang.
\newblock {\em Projected dynamical systems and variational inequalities with
  applications}, volume~2.
\newblock Springer Science \& Business Media, 2012.

\bibitem{NocWri00}
J.~Nocedal and S.~Wright.
\newblock {\em Numerical Optimization}.
\newblock Springer Series in Operation Research and Financial Engineering. 2nd
  edition, 2000.

\bibitem{rockafellar2}
R.~T. Rockafellar and R.~J.-B. Wets.
\newblock {\em Variational Analysis}.
\newblock Springer-Verlag Berlin Heidelberg, 1998.

\bibitem{Seron1997}
M.M. Seron, J.H. Braslavsky, and G.C. Goodwin.
\newblock {\em {Fundamental Limitations in Filtering and Control}}.
\newblock Springer London, London, 1997.

\bibitem{ShaHee_ACC21a}
B.~Sharif, M.~Heertjes, H.~Nijmeijer, and W.P.M.H. Heemels.
\newblock On the equivalence of extended and oblique projected dynamics with
  applications to hybrid integrator-gain systems.
\newblock In {\em American Control Conference (ACC) 2021, New Orleans, USA},
  pages 3434--3439, 2021.

\bibitem{ShaHee_CDC19a}
B.~Sharif, M.F. Heertjes, and W.P.M.H. Heemels.
\newblock Extended projected dynamical systems with applications to hybrid
  integrator-gain system.
\newblock In {\em IEEE Conf.~Decision and Control (CDC)}, pages 5773--5778,
  2019.

\bibitem{ShiNi_CDC22}
K.~Shi, N.~Nikooienejad, I.R. Petersen, and S.~O.~R. Moheimani.
\newblock A negative imaginary approach to hybrid integrator-gain system
  control.
\newblock In {\em IEEE Conference Decision and Control}, pages 1968--1973,
  2022.

\bibitem{EijHee_CSL20a}
S.J.A.M. van~den Eijnden, M.F. Heertjes, W.P.M.H. Heemels, and H.~Nijmeijer.
\newblock Hybrid integrator-gain systems: A remedy for overshoot limitations in
  linear control?
\newblock {\em L-CSS}, 4:1042 -- 1047, 2020.

\end{thebibliography}

\end{document}